\documentclass{amsart}

\usepackage {amsfonts}
\usepackage[english]{babel}
\def\g{\gamma}
\def\G{\Gamma}
\def\d{\delta}
\def\a{\alpha}
\def\b{\beta}
\def\p{\varphi}
\def\e{\varepsilon}
\def\l{\lambda}
\def\L{\Lambda}
\def\s{\sigma}

\def\o{\omega}
\def\O{\Omega}

\def\R{{\mathbb R}}
\def\C{{\mathbb C}}
\def\N{{\mathbb N}}
\def\Z{{\mathbb Z}}

\def\Im{\mbox{Im }}
\def\bs{~\hfill\rule{7pt}{7pt}}
\def\F{\frak T}

\DeclareMathOperator{\spec}{sp}
\DeclareMathOperator{\Res}{Res}
\DeclareMathOperator{\const}{const}
\DeclareMathOperator{\supp}{supp }
\DeclareMathOperator{\dist}{dist}

\newtheorem{Th}{Theorem}
\newtheorem{Pro}{Proposition}
\newtheorem{Def}{Definition}

\newtheorem{Cor}{Corollary}
\begin{document}

\title{Generalized Fourier quasicrystals, almost periodic sets, and zeros of Dirichlet series}

\author{Sergii Yu.Favorov}

\address{Sergii Favorov,
\newline\hphantom{iii}  Karazin's Kharkiv National University
\newline\hphantom{iii} Svobody sq., 4, Kharkiv, Ukraine 61022}
\email{sfavorov@gmail.com}

\maketitle {\small
\begin{quote}
\noindent{\bf Abstract.}
Let $S(z)$ be an absolutely convergent Dirichlet series with a bounded spectrum and only real zeros $a_n,\,n=1,2,\dots$, let $\mu$ be the sum of unit masses at points $a_n$. It is proven that the Fourier transform $\hat\mu$ in the sense of distributions is a purely point measure. Conversely, in terms of the properties of $\hat\mu$, a sufficient condition was found when $a_n$ are zeros of an absolutely convergent Dirichlet series with bounded spectrum. Also, in terms of the properties of $\hat\mu$, a criterion is established that $a_n$ are zeros of an almost periodic entire function of exponential growth. In all cases, the multiplicity of zeros is taken into account.

  Almost periodic sets, introduced by B.Levin and M.Krein in 1948, play important role in our investigations. In particular, we show a new simple representation of such sets.

\medskip


\end{quote}
}

\medskip

   \section{Introduction}\label{S1}
   \bigskip

A crystalline measure on $\R^d$ is a complex measure  with discrete locally finite support, which is a temperate distribution and its distributional Fourier transform $\hat\mu$
is also a measure with locally finite support; if, in addition, the measures $|\mu|$ and $|\hat\mu|$ are temperate distributions, then $\mu$ is called the Fourier quasicrystal.

The  Fourier quasicrystal may be considered as a mathematical model for  atomic arrangement
having a discrete diffraction pattern.  There are a lot of papers devoted to study properties of Fourier quasicrystals or, more generally, crystalline measures.
For example, one can mark collections of papers \cite{D}, \cite{Q}, in particular, the basic paper \cite{L1}.

Measures of the form
\begin{equation}\label{a}
\mu=\sum_{\l\in\L} c_\l\d_\l,\qquad c_\l\in\N,
\end{equation}
  are the most important case of Fourier quasicrystal. Recently A. Olevsky and A. Ulanovsky \cite{OU1}, \cite{OU} showed that the supports $\L$ of these measures are zero sets of exponential polynomials with purely imaginary exponents and only real zeros with multiplicities $c_\l$ at points $\l\in\L$ and, conversely, the zero sets of such exponential polynomials are supports of some Fourier quasicrystals of the form \eqref{a}.

In our article, using the same methods, we present analogues results for measures \eqref{a} with the distribution Fourier transform
\begin{equation}\label{b}
  \hat\mu=\sum_{\g\in\G}b_\g\d_\g,
\end{equation}
where $\G$ is an arbitrary countable set. In this case the corresponding Poisson's formula
$$
\sum_{\l\in\L} c_\l\hat f(\l)=\sum_{\g\in\G}b_\g f(\g).
$$
also takes place for every function $f$ from Schwarz's class. In order to describe such measures, we use the concept of almost periodic sets, which was introduced
 by M. Krein and B. Levin (\cite{L}, App. VI). In modern notations (cf.\cite{M1}, \cite{R}), a locally finite set $\L$ with multiplicities $c_\l$ at points $\l\in\L$
 is almost periodic, if the convolution of measure \eqref{a} with every continuous function with compact support is an almost periodic function.
 We will write an almost periodic set like the sequence $A=\{a_n\}_{n\in\Z}$, where each point $a_n=\l$ occurs $c_\l$ times. Therefore, almost periodic sets are in fact multisets.

In Section \ref{S2} we give the original definition of almost periodic sets given by Krein and Levin, which is equivalent to the one given above.
Also we prove some properties of almost periodic sets, in particular, show that such  sets
have the form $\{\a n+\phi(n)\}_{n\in\N}$ with $\a>0$ and an almost periodic mapping $\phi:\,\Z\to\R$.

In Section \ref{S3} we consider an absolutely convergent Dirichlet series with bounded spectrum and a real set of zeros. We prove that the Fourier transform of the measure $\mu_A$ corresponding to such a set $A$ is always a purely point measure. Note that the zero set of any absolutely convergent Dirichlet series (or, more general, any almost periodic holomorphic function and even any holomorphic function with almost periodic modulus) is almost periodic (cf.\cite{FRR}).

In section \ref{S4} we investigate the inverse problem. We consider a measure $\mu$ of form \eqref{a} with the Fourier transform $\hat\mu$ of form \eqref{b} and find sufficient conditions
for which $A$ is the zero set of an absolutely convergent Dirichlet series with bounded spectrum.

Here we assume that the measure $|\hat\mu|$ is a temperate distribution. By \cite {F4}, Lemma 1, the multiset $A=\{a_n\}_{n\in\Z}$, where each point $a_n=\l$ occurs $c_\l$ times, is  almost periodic.
 In \cite{FRR} it was proven that every almost periodic set $A\subset\R$ is exactly the zero set of some entire almost periodic function. Every almost periodic function is bounded on the real line,
hence, by the Phragmen-Lindelöf principle, any almost periodic entire function cannot grow less than exponentially. We  find a criterion for  $A$ to be the zero set
of an almost periodic entire function of the exponential growth in terms of $\g$ and $\b_\g$ from \eqref{b}.

\bigskip

\section{Almost periodic sets}\label{S2}
\bigskip

\begin{Def}[for example, see \cite{B}, \cite{Le}]
 A continuous function $g(x)$ on the real line
is  almost periodic if for any  $\e>0$ the set of $\e$-almost periods
  $$
E_{\e}= \{\tau\in\R:\,\sup_{x\in\R}|g(x+\tau)-g(x)|<\e\}
  $$
is relatively dense, i.e., $E_{\e}\cap(x,x+L)\neq\emptyset$ for all $x\in\R$ and some $L$ depending on $\e$.
\end{Def}
For example every sum
$$
Q(x)=\sum q_n e^{2\pi ix\o_n},\quad\o_n\in\R,\quad q_n\in\C,\qquad \sum_n |q_n|<\infty,
$$
is an almost periodic function.

Spectrum  of an almost periodic function $g$ is the set
  $$
   \spec g=\left\{w\in\R:\,\lim_{t\to\infty}\frac{1}{2T}\int_{-T}^Te^{-2\pi i\o x} g(x)dx\neq0\right\}.
  $$
It is easy to see that $\spec Q(x)=\{\o_n:\,q_n\neq0\}$. Note that spectra of almost periodic functions are at most countable.

\begin{Def}[see \cite{B}, \cite{Le}]
 A continuous function $g(z)$ on a strip
$$
S=\{z=x+iy:\,-\infty\le a<y<b\le+\infty\}\subset\C
$$
is  almost periodic if for any $\a,\,\b$ such that $[\a,\b]\subset(a,b)$ and
  $\e>0$ the set of $\e$-almost periods
  $$
E_{\a,\b,\e}= \{\tau\in\R:\,\sup_{x\in\R,\a\le y\le\b}|g(x+\tau+iy)-g(x+iy)|<\e\}
  $$
is relatively dense, i.e., $E_{\a,\b,\e}\cap(x,x+L)\neq\emptyset$ for all $x\in\R$ and some $L$ depending on $\e,\a,\b$.
\end{Def}

Just as it was done in \cite{L}, we could define an almost periodic set in a strip. Here we give
this definition only in the simplified form. The general definition, which takes into account the behavior of $A$ near the boundary, was given by H.Tornehave in \cite{T}.
The connection between almost periodic sets of general form and zeros of holomorphic almost periodic functions in terms of Chern cohomology see \cite{F}.

\begin{Def}[M. Krein and B. Levin \cite{L}, App. VI]\label{D2}
 Let $S$ be a horizontal strip of a finite width. A discrete locally finite multiset $A=\{a_n\}_{n\in\Z}\subset S$ is almost periodic if for any $\e>0$ the set of its $\e$-almost periods
\begin{equation}\label{e}
E_\e=\{\tau\in\R:\,\exists\ \text{ a bijection } \s:\Z\to\Z \quad\text{such that}\ \sup_n |a_n+\tau-a_{\s(n)}|<\e\}
\end{equation}
 has nonempty intersection with every interval $(x,x+L_\e)$.
\end{Def}
In our article we consider, as a rule, only the case of almost periodic sets on the real line.

\medskip
Set $\mu_A=\sum_n \d_{a_n}$. Clearly, the mass of $\mu_A$ in any point $x\in\R$ is equal to the multiplicity of this point in the sequence $\{a_n\}_{n\in\Z}$.

It was proved in \cite{FRR} that almost periodicity of $A$ is equivalent to almost periodicity of convolution $\mu_A\star\p$ for every $C^\infty$-function $\p(x),\,x\in\R$,
 with compact support. But it is easy to replace $C^\infty$-functions by continuous functions with compact support.

 Indeed, take $C^\infty$-function $\p\ge0$ such that $\p(x)\equiv1$ for $0<x<1$. If $\mu_A\star\p$ is almost periodic, then it is uniformly bounded, hence
 $\mu_A[x,x+1]<k_1$ for all $x\in\R$ with some constant $k_1$. For any continuous function $\psi$ with support in $(0,1)$ one can take $\p\in C^\infty$ such that
 $\sup_{x\in\R}|\psi(x)-\p(x)|<\e/k_1$. We obtain that every $\e$-almost period of $\mu_A\star\p$ is $2\e$-almost period of $\mu_A\star\psi$.

By the way, we gave the proof of the following proposition
\begin{Pro}[\cite{FRR}]\label{P1}
For any almost periodic set there is $k_1\in\N$ such that $\# A\cap[x,x+1]\le k_1$. Also, $\# A\cap[x,x+h)\le k_1(h+1)$.
\end{Pro}
Here and below, $\# H$ means the number of points in the multiset $H$ taking into account their multiplicities.
\smallskip

\begin{Pro}\label{P2}
For any almost periodic set there is $k_2\in\N$ such that for every $h>0$ and every half-intervals $[x_1,x_1+h),\,[x_2,x_2+h)$ we have
 $|\# A\cap[x_1,x_1+h)-\# A\cap[x_2,x_2+h)|\le k_2$. Also, for every $x\in\R,\, h>0,\,M\in\N$
$$
|\# A\cap[x,x+h)-(1/M)\# A\cap[x,x+Mh)|\le k_2.
$$
\end{Pro}
{\bf Proof}. Let $L_1,\,E_1$ be defined in \eqref{e}, and $\tau\in E_1\cap[x_1-x_2,L_1+x_1-x_2)$. Since
$[x_2,x_2+h)+\tau\subset[x_1,x_1+L_1+h)$, we see that to each $a_n\in[x_2,x_2+h)$ assign a point $a_{\s(n)}\in[x_1-1,x_1+L_1+h+1)$. Therefore,
$$
\# A\cap[x_2,x_2+h)\le \# A\cap[x_1,x_1+h)+\# A\cap[x_1-1,x_1)+\# A\cap[x_1+h,x_1+h+L_1+1).
$$
By Proposition \ref{P1}, the last two terms are bounded by $k_1+(L_1+2)k_1$. The proof of the opposite inequality is the same.

To prove the second assertion, we have to add together all the inequalities
$$
\# A\cap[x,x+h)-k_2\le\# A\cap[x+(m-1)h,x+mh)\le k_2+\# A\cap[x,x+h),\quad m=1,2,\dots,M.
$$
\bs

\begin{Pro}\label{P3}
Let $A$ be an almost periodic set. There is a strictly positive density $d$ such that for any $\eta>0$ and any half-interval $I$  with length $l(I)>N_\eta$ we have
$$
            \left|\frac{\# A\cap I}{l(I)}-d\right|<\eta.
$$
\end{Pro}
This result was generalized to multidimensional Euclidean spaces in \cite{FK2}.
\medskip

{\bf Proof of Proposition \ref{P3}}. Let $I_1=[x_1,x_1+h_1),I_2=[x_2,x_2+h_2)$ be two half-intervals such that $h_1/h_2=p/q,\,p,q\in\N$. We have
$$
  \frac{\# A\cap I_1}{h_1}-\frac{\# A\cap I_2}{h_2}=\frac{\# A\cap I_1}{h_1}-\frac{\# A\cap qI_1}{qh_1}+\frac{\# A\cap qI_1}{qh_1}-\frac{\# A\cap pI_2}{ph_2}+
  \frac{\# A\cap pI_2}{ph_2}-\frac{\# A\cap I_2}{h_2}.
$$
Applying Proposition \ref{P2}, we get
\begin{equation}\label{i}
 \left|\frac{\# A\cap I_1}{h_1}-\frac{\# A\cap I_2}{h_2}\right|\le \frac{k_2}{h_1}+\frac{k_2}{qh_1}+\frac{k_2}{h_2}\le k_2\left(\frac{2}{h_1}+\frac{1}{h_2}\right).
\end{equation}
    For arbitrary $h_1,\,h_2$ take a half-interval $I'=[x_1,x_1+h')$ such that $h_1<h'<h_1+1$ and $h'/h_2$ rational. We have
$$
 \left|\frac{\# A\cap I_1}{h_1}-\frac{\# A\cap I'}{h'}\right|\le\frac{\# A\cap[x_1+h_1,x_1+h')}{h'}+ \frac{\#A\cap I_1}{h_1h'}.
$$
By Proposition \ref{P1}, we obtain
$$
 \left|\frac{\# A\cap I_1}{h_1}-\frac{\# A\cap I'}{h'}\right|\le\frac{k_1}{h'}+\frac{k_1(h_1+1)}{h_1h'}.
$$
Applying \eqref{i} with $I'$ instead of $I_1$, we obtain for all $I_1,\,I_2$
$$
 \left|\frac{\# A\cap I_1}{l(I_1)}-\frac{\# A\cap I_2}{l(I_2)}\right|\le k_2\left(\frac{2}{l(I_1)}+\frac{1}{l(I_2)}\right)+k_1\left(\frac{2}{l(I_1)}+\frac{1}{l(I_1)^2}\right).
$$
Therefore there is a limit
$$
  d=\lim_{l(I)\to\infty}\frac{\# A\cap I}{l(I)}.
$$
It is easy to check that the set $A$ has nonempty intersection with every interval of length $2+L_1$, hence this limit is strictly positive.  \bs
\medskip

\begin{Th}\label{T1}
Let $A=\{a_n\}\subset\R$ be an almost periodic set of density $d$ such that $a_n\le a_{n+1}$ for all $n\in\Z$.
Then
\begin{equation}\label{p}
a_n=n/d+\phi(n)\quad\text{with an almost periodic mapping}\quad \phi:\,Z\to\R.
\end{equation}
\end{Th}
{\bf Remark 1}. The incomplete proof of this Theorem was given in \cite{FK1}.
\medskip

{\bf Remark 2}. The converse assertion is simple, because for every $\e$-almost period $\tau\in\Z$ of the mapping $\phi(n)$ the number $\tau/d$ is $\e$-almost period of an almost periodic set $A=\{n/d+\phi(n)\}$
with any $d>0$ and $\s(n)=n+\tau$.
\medskip

{\bf Proof of Theorem \ref{T1}}. We may suppose that $a_0<a_1$.
It follows from Proposition \ref{P1} that every interval of length $1$ contains at least one subinterval of length $1/(2k_1)$ that does not intersect $A$.
Take $\e<\min\{1/(6k_1),(a_1-a_0)/3\}$. Divide $\R$ into an infinite number of disjoint half-intervals $I_j=(t_j,t_{j+1}],\,j\in\Z$ such that $t_{j+1}-t_j<2$
and $A\cap (t_j-2\e,t_j+2\e)=\emptyset$ for all $j$.

 Let $\tau$ be any number from $E_\e$ in \eqref{e}, and let $\s$ be the corresponding  bijection.  Then $\rho(j)\in\Z$ corresponds to any $j$
such that $\s$ is the bijection of $A\cap I_j$ to $A\cap I_{\rho(j)}$. Hence, $\#(A\cap I_j)=\#(A\cap I_{\rho(j)})$. Let $\s_j$ be the monotone increasing bijection of $A\cap I_j$ on $A\cap I_{\rho(j)}$. Check that
\begin{equation}\label{so}
|a_n+\tau-a_{\s_j(n)}|<\e\qquad\forall\, a_n\in I_j.
\end{equation}
   Suppose the contrary.  Let $n_0$ be the minimal number such that \eqref{so} does not satisfy. If $a_{n_0}+\tau+\e\le a_{\s_j(n_0)}$, then
$a_n+\tau+\e\le a_k$ for all $n\le n_0$ and $k\ge\s_j(n_0)$, $a_n\in I_j,\,a_k\in I_{\rho(j)}$. Therefore, $k\neq \s(n)$ for these numbers,
and $\s$ gives  a correspondence
between points from the set $\{n\le n_0:\,a_n\in I_j\}$ and points from the set $\{k<\s_j(n_0):\,a_k\in I_{\rho(j)}\}$. But
by definition of $\s_j$, we have
$$
\#\{n\le n_0:\,a_n\in I_j\}=\#\{k\le\s_j(n_0):\,a_k\in I_{\rho(j)}\}=1+\#\{k<\s_j(n_0):\,a_k\in I_{\rho(j)}\}.
$$
We get a contradiction.

If $a_{n_0}+\tau\ge a_{\s_j(n_0)}+\e$, then $a_n+\tau\ge a_k+\e$ for all $n\ge n_0$ and $k\le \s_j(n_0)$, $a_n\in I_j,\,a_k\in I_{\rho(j)}$.
Therefore, $k\neq \s(n)$ for these numbers, and $\s$ gives  a correspondence
between points from the set $\{n\ge n_0:\,a_n\in I_j\}$ and points from the set $\{k>\s_j(n_0):\,a_k\in I_{\rho(j)}\}$. But by definition of $\s_j$, we have
$$
\#\{n\ge n_0:\,a_n\in I_j\}=\#\{k\ge\s_j(n_0):\,a_k\in I_{\rho(j)}\}=1+\#\{k>\s_j(n_0):\,a_k\in I_{\rho(j)}\}.
$$
We get a contradiction as well.

Since numbers $\#(A\cap I_j)$ and $\#(A\cap I_{\rho(j)})$ coincide, we see that the differences between indices of the first elements in these sets coincide for all $j$.
Hence there is a number $h\in\Z$ so that  inequality \eqref{e} satisfies for all $n\in\N$ with $\s(n)=n+h$.

It follows from the definition of $\tau$ for all $k\in\N$ and $N\in\N$
$$
   -\e<a_{kh}-\tau-a_{(k-1)h}<\e, \quad\text{and}\quad  -N\e<a_{Nh}-N\tau-a_0<N\e.
$$
Let $I$ be the smallest half-interval containing $a_0$ and $a_{Nh}$. The last inequality implies that its length satisfies the inequality
$$
 N\tau-N\e<l(I)<N\tau+N\e.
$$
 On the other hand, taking into account that ends of $I$ may be points of $A$ with multiplicity at most $k_1$, we have
$$
 Nh-2(k_1-1)\le\#A\cap I\le Nh+2(k_1-1).
$$
Therefore,
$$
\frac{Nh-2(k_1-1)}{N\tau+N\e}\le\frac{\#A\cap I}{l(I)}\le\frac{Nh+2(k_1-1)}{N\tau-N\e}.
$$
Passing to the limit as $N\to\infty$ and using Proposition \ref{P3}, we obtain  the inequality
$$
  \tau-\e\le h/d\le\tau+\e.
$$
 Set $\phi(n):=a_n-n/d$. We get for all $n\in\Z$
$$
  \phi(n+h)-\phi(n)=a_{n+h}-(n+h)/d-a_n+n/d=a_{\s(n)}-(a_n+\tau)+(\tau-h/d).
$$
Using \eqref{e}, we obtain $|\phi(n+h)-\phi(n)|<2\e$. Therefore, $h$ is $2\e$-almost period of the function $\phi$. The set of $\e$-almost periods $\tau$ of $A$
 is relatively dense, therefore the set of such integers $h$ is relatively dense as well.  \bs

\begin{Cor}\label{C1}
For any almost periodic set $A=\{a_n\}$ such that $0\not\in A$ there is a finite limit
$$
\lim_{N\to\infty}\sum_{|a_n|<N}1/a_n.
$$
Moreover,  the sum
$$
\frac{1}{z-a_0}+\sum_{n\in\N\setminus\{0\}}\left[\frac{1}{z-a_n}+\frac{1}{z-a_{-n}}\right].
$$
converges absolutely and uniformly on every disjoint with $A$ compact set $K$.
\end{Cor}
{\bf Proof}. Let $A=\{n/d+\phi(n)\}_{n\in\Z}$. Since the numbers $\phi(n)$ are uniformly bounded, we see that  the sums
$$
 \sum_{n\in\Z,|a_n|<N}\frac{1}{a_n} \quad\mbox{and}\quad \sum_{n\in\Z,|n|<dN}\frac{1}{n/d+\phi(n)}
$$
 differ for a uniformly bounded with respect to $N$ number of terms, and each of these terms tends to $0$ as $N\to\infty$. Then
 $$
   \sum_{n\in\Z,0<|n|<N}\frac{1}{n/d+\phi(n)}=\sum_{n\in\N,0<n<N}\frac{\phi(n)+\phi(-n)}{\phi(n)\phi(-n)+n\phi(-n)/d-n\phi(n)/d-(n/d)^2}.
 $$
 The first assertion follows from Cauchy criterium. The second one follows from the absolutely convergence of the series
 $$
 \sum_{n\in\N\setminus\{0\}}\left[\frac{1}{z-a_n}+\frac{1}{z-a_{-n}}\right]=\sum_{n\in\N\setminus\{0\}}\left[\frac{2z-\phi(-n)-\phi(n)}{(n/d+\phi(n)-z)(-n/d+\phi(-n)-z)}\right].
 \phantom{XXX} \bs
 $$
 \medskip

In \cite{L}, App. VI, M.Krein and B.Levin considered zero sets $Z_f$ of entire almost periodic functions $f$ of exponential growth. They proved that if
$Z_f\subset\R$, then its zeros $a_n$ form an almost periodic set, which satisfy \eqref{p} and
\begin{equation}\label{p1}
 \sup_{\tau\in\Z}\sum_{n\in\Z\setminus\{0\}}n^{-1}[\phi(n+\tau)-\phi(n)]<\infty.
\end{equation}
On the other hand, they proved that any almost periodic set $A\subset\R$ satisfying conditions \eqref{p} and \eqref{p1} is the set of zeros
of an entire almost periodic function of exponential growth.

It follows from Theorem \ref{T1} that condition \eqref{p} can be omitted in the last result.\medskip

Theorem \ref{T1} was generalized by W.Lawton \cite{La} to almost periodic sets in $\R^m,\,m>1,$ whose spectrum is contained
 in a finitely generated  additive group.
\bigskip

 \section{Zeros of infinite exponential sums}\label{S3}
\bigskip

Denote by $S(\R)$ the Schwartz space of test functions $\p\in C^\infty(\R)$ with the finite norms
 $$
  N_{n,m}(\p)=\sup_{\R}\max_{k\le m} |(1+|x|^n)\p^{(k)}(x)|,\quad n,m=0,1,2,\dots
 $$
 These norms generate the topology on $S(\R)$.  Elements of the space $S^*(\R)$ of continuous linear functionals on $S(\R)$ are called temperate distributions.

The Fourier transform of a temperate distribution $f$ is defined by the equality
$$
\hat f(\p)=f(\hat\p)\quad\mbox{for all}\quad\p\in S(\R),
$$
where
$$
   \hat\p(t)=\int_{\R^d}\p(x)e^{-2\pi i xt}dx
 $$
is the Fourier transform of the function $\p$. By $\check\p$ we  denote the inverse Fourier transform of $\p$.
The Fourier transform is the bijection of $S(\R)$ on itself and the bijection of $S^*(\R)$ on itself.

\medskip
Let $\F$ be the class of exponential sums
$$
f(x)=\sum_n q_n e^{2\pi i\o_n x},\quad q_j\in\C\setminus\{0\}
$$
 with finite Wiener's norm $\|f\|_W:=\sum_n|q_n|$ and a bounded spectrum $\O:=\{\o_n\}\subset\R$.

 Every function $f(x)\in\F$ expands to an entire almost periodic function $f(z)$ of exponential type $\s=\max_n|\o_n|$; the zero set $A=\{a_n\}$ of $f(z)$
 lies in some horizontal strip of a finite width if and only if $\inf\O\in\O$ and $\sup\O\in\O$ (cf.\cite{L}, Ch.VI, Cor.2). Moreover, $A$ is an almost periodic set (cf.\cite{L}, App.VI, Lemma 1).

 If $0\not\in A$ and  $a_n=\a n+\phi(n)$ with almost periodic $\phi:\,\Z\to\R$ such that the Fourier series of $\phi$ converges absolutely, then the function
$$
  (1-z/a_0)\prod_{n\in\N} (1-z/a_n)(1-z/a_{-n})
$$
expands in an absolutely convergent exponential series (cf.\cite{L}, App. VI, Th.9), hence it belongs to $\F$.
\begin{Th}\label{T2}
Suppose that $f\in\F$ has the zero set $A=\{a_n\}\subset\R$, and $\mu_A=\sum_n\d_{a_n}$. Then the Fourier transform $\hat\mu_A$ is a pure point measure.
\end{Th}
{\bf Remark}. Y.Meyer (\cite{M}, Th.5.16) proved that under condition  $a_n=\a n+\phi(n)$ with an absolutely convergent Fourier series for $\phi$ the temperate distribution $\hat\mu_A$ is a pure point measure.
It follows from \cite{L}, App. VI, Th.9, that the above Theorem \ref{T2} is a strengthening of Meyer's result.
\medskip

 {\bf Proof}. It follows from Proposition \ref{P1} that the measure $\mu_A$ satisfies the condition $\mu_A([-r,r])=O(r)$ as $r\to\infty$. Hence the measure $\mu_A$ and  the distribution $\hat\mu_A$ are temperate
 distributions.
 In order to prove that  $\hat\mu_A$ is a measure, we will check the estimate
 \begin{equation}\label{b1}
 |(\hat\mu_A,\p)|\le C\max_{|t|<T}|\p(t)|
 \end{equation}
  for any $T<\infty$ and any $C^\infty$-function $\p$ with support on the interval $(-T,T)$. If this is the case, the distribution $\hat\mu$ has a unique expansion to a linear functional on the space of continuous functions $g$ on $[-T,T]$ such that $g(\a)=g(\b)=0$  with bound \eqref{b1}. Since we can expand this functional to the space of all continuous functions on $[-T,T]$ with bound \eqref{b1} too, we see that $\hat\mu_A$ is a complex measure.

Let $\p$ be $C^\infty$-function with support in $(-T,T)$. Set
$$
\Phi(z)=\int_{-\infty}^\infty \p(t)e^{-2\pi itz}dt.
$$
Clearly,  $\Phi(z)$ is an entire function, which equals the Fourier transform of the function $\p(t)e^{2\pi ty}$.
Therefore, $\Phi(x+iy)$  belongs to $S(\R)$ for each fixed $y\in\R$, and  we have for its converse Fourier transform
\begin{equation}\label{cF}
 \check\Phi(\o+iy)=\p(\o)e^{2\pi\o y},\qquad \o\in\R.
\end{equation}
Then $f(z)$ is almost periodic,  hence Lemma 1 from \cite{L}, Ch.6, implies that for every $\e>0$ and $s<\infty$ there is $m=m(\e,s)>0$ such that
 $$
 |f(z)|\ge m\quad\mbox{for}\quad |\Im z|\le s\quad\mbox{and}\quad z\not\in A(\e):=\{z:\,\dist(z,A)<\e\}.
 $$
  By Proposition \ref{P1}, for $\e$ small enough each connected component of $A(\e)$ contains no segment of length $1$, hence its diameter is less than $1$.
Therefore there are two sequences
$R_k\to+\infty,\,R'_k\to-\infty$ such that
$$
   |f(x+iy)|>m\quad\text{for}\quad x=R_k \quad\text{or}\quad x=R'_k,\quad |y|<1.
$$
 Consider integrals of the function $\Phi(z)f'(z)f^{-1}(z)$ over the boundaries of the rectangles
$\Pi_k=\{z=x+iy:\,R'_k<x<R_k,-s'<y<s\}$, where the numbers $s,\,s'>0$ will be chosen later. Since $\Phi(x\pm iy)$ tends to zero as $x=R_k\to+\infty,\,x=R'_k\to-\infty$ uniformly with respect to $-s'\le y\le s$, we obtain that these integrals tend to
\begin{equation}\label{in}
   \int_{+\infty}^{-\infty}\Phi(x+is)f'(x+is)f^{-1}(x+is)dx+\int_{-\infty}^{+\infty}\Phi(x-is')f'(x-is')f^{-1}(x-is')dx=:I_1+I_2..
\end{equation}
 By the Theorem on residues
 \begin{equation}\label{r}
 I_1+I_2= 2\pi i\sum_{\l:f(\l)=0}\Res_\l \Phi(z)f'(z)f^{-1}(z)=2\pi i\sum_{\l:f(\l)=0}a(\l)\Phi(\l)=2\pi i(\mu_A,\Phi),
\end{equation}
where $a(\l)$ is the multiplicity of $f(z)$ at the point $\l$.

Set $\o_1=\inf\O$. Then the corresponding coefficient $q_1$ does not vanish. Taking into account that $\sum_n|q_n|<\infty$, we can take numbers  first $M$ and then $s$ such that
\begin{equation}\label{es}
\sum_{n>M}|q_n/q_1|<1/3,\phantom{XXXXXXXX} e^{2\pi(\o_1-\o')s}\sum_{n=2}^\infty |q_n/q_1|<1/3,
\end{equation}
where $\o'=\min_{2\le n\le M}\o_n$. We have  for $z=x+is$
$$
   f(z)=q_1 e^{2\pi i(x+is)\o_1}\left(1+\sum_{n=2}^M (q_n/q_1)e^{2\pi i(\o_n-\o_1)(x+is)}+\sum_{n>M} (q_n/q_1)e^{2\pi i(\o_n-\o_1)(x+is)}\right).
$$
   Set
$$
  H(x):= \sum_{n=2}^\infty (q_n/q_1)e^{2\pi i(\o_n-\o_1)x}e^{2\pi(\o_1-\o_n)s}=\sum_{n=2}^\infty h_n(s)e^{2\pi i(\o_n-\o_1)x}.
$$
 By \eqref{es}, we have $\|H\|_W<2/3$. Since $\|\cdot\|_W$ is the norm on the algebra of all absolutely convergent exponential sums, we get
\begin{equation}\label{H}
(1+H(x))^{-1}=\sum_{j=0}^\infty(-1)^j H^j(x),\qquad \|(1+H)^{-1}\|_W\le\sum_{j=0}^\infty\|H^j\|_W<3.
\end{equation}

We have
$$
f^{-1}(x+is)=q_1^{-1}e^{2\pi\o_1 s}e^{-2\pi i\o_1x}(1+H(x))^{-1},\qquad f'(x+is)=\sum_{n=1}^\infty 2\pi i\o_nq_n e^{-2\pi\o_ns}e^{2\pi i\o_nx},
$$
and
\begin{equation}\label{H1}
  f'(x+is)f^{-1}(x+is)=\sum_{n=1}^\infty 2\pi i\o_n(q_n/q_1) e^{-2\pi(\o_n-\o_1)s}e^{2\pi i(\o_n-\o_1)x}(1+H(x))^{-1}.
\end{equation}
Rewrite $f' f^{-1}$  in the form
$$
 f'(x+is)f^{-1}(x+is)=\sum_{\g\in\G_1} p_\g e^{2\pi i\g x},\quad p_\g=p_\g(s)\in\C,
$$
 with some countable $\G_1\subset\R_+\cup\{0\}$. Since $\O$ is bounded, we obtain from \eqref{H} and \eqref{H1}
$$
 \sum_{\g\in\G_1}|p_\g|=\|f^{-1}(x+is) f'(x+is)\|_W\le 6\pi\max_n\{|\o_n| e^{2\pi(\o_n-\o_1)s}\} \sum_n|q_n/q_1| =:C_f.
$$
  The function $\Phi(x+is)$ belongs to $S(\R)$ for $s$ fixed, hence $|x|^2\Phi(x+is)\to0$ as $|x|\to\infty$.
 Changing the order of integration and summation and taking into account \eqref{cF}, we obtain for  the first integral in \eqref{in}
 $$
  I_1= -\sum_{\g\in\G_1}p_\g\int_{-\infty}^{+\infty}\Phi(x+is)e^{2\pi i\g x}dx=-\sum_{\g\in\G_1}p_\g e^{2\pi\g s}\p(\g).
$$
Since $\supp\p\subset(-T,T)$, we get the bound
$$
  |I_1|\le C_f e^{2\pi Ts}\max_{|t|\le T}|\p(t)|.
$$
The similar arguments show that the second integral in \eqref{in} with the appropriate $s'$ has the same bound.

Since $\hat\p(x)=\Phi(x)$, we get from \eqref{r}
$$
(\mu_A,\hat\p)=(2\pi i)^{-1}(I_1+I_2).
$$
Therefore,
$$
|(\hat\mu_A,\p)|=|(\mu_A,\hat\p)|\le C(f,T)\sup_{|y|\le T}|\p(y)|,
$$
and $\hat\mu_A$ is a measure. Since $\mu_A$ is almost periodic, Theorem 5.5 from \cite{M} implies that $\hat\mu_A$ is a pure point measure.
\bs

 \section{Entire functions with the given almost periodic zero sets}\label{S4}
\bigskip

In this section we assume that  a measure $\mu$ of form \eqref{a} is a temperate distribution, its Fourier transform $\hat\mu$ is a pure point measure of form \eqref{b}, the measure $|\hat\mu|$ is  a temperate distribution, and
$A=\{a_n\}_{n\in\Z}$ is a multiset, where each point $a_n=\l\in\supp\mu$ occurs $c_\l$ times. By \cite {F4}, Lemma 1,  $\mu$ is an almost periodic measure and $A$ is an almost periodic set. Since $\hat\mu$ is also a measure,
Theorem 5.5 from \cite{M} implies that every number $b_\g$ from \eqref{b} equals the corresponding Fourier coefficient of the measure $\mu$, i.e.,
$$
  b_\g=\lim_{T\to\infty}\frac{1}{2T}\int_{-T}^T e^{-2\pi i\g x} \mu(dx).
$$
In particular, $b_{-\g}=\bar b_\g$ and $b_0\ge0$. Moreover,  $b_0$ coincides with the density $d$ of the set $A$. Also, it can be checked (cf.\cite{F1}) that the condition $|\hat\mu|\in S^*(\R)$ implies
 \begin{equation}\label{ch}
|\hat\mu|(-r,r)= \sum_{|\g|<r}|b_\g|=O(r^\kappa)\quad\text{as}\quad r\to\infty
\end{equation}
with some $\kappa<\infty$.

In what follows we will suppose that $0\not\in A$. By Corollary \ref{C1}, the set $A$ satisfies Lindelof's condition.  Hence, the function
\begin{equation}\label{f}
  f(z)=(1-z/a_0)\prod_{n\in\N} (1-z/a_n)(1-z/a_{-n})
\end{equation}
 is an entire function of exponential type with the zero set $A$. Note that $f(\bar z)=\overline{f(z)}$.
Also, introduce the notations
 $$
 \R_+:=\{x\in\R:x>0\},\,\R_-:=-\R_+,\,\C_+:=\{z\in\C:\Im z>0\},\,\C_-:=-\C_+.
 $$
\begin{Pro}\label{P4}
For all $z=x+iy\in\C_+$
\begin{equation}\label{n1}
\frac{f'(z)}{f(z)}=\frac{1}{z-a_0}+\sum_{n\in\N}\left[\frac{1}{z-a_n}+\frac{1}{z-a_{-n}}\right]=-2\pi i\sum_{\g\in\G\cap\R_+}b_\g e^{2\pi i\g z}-\pi id,
\end{equation}
and for all $z=x+iy\in\C_-$
\begin{equation}\label{n2}
\frac{f'(z)}{f(z)}=\frac{1}{z-a_0}+\sum_{n\in\N}\left[\frac{1}{z-a_n}+\frac{1}{z-a_{-n}}\right]=2\pi i\sum_{\g\in\G\cap\R_-}b_\g e^{2\pi i\g z}+\pi id,
\end{equation}
where $d$ is the density of the almost periodic set $A$.

Then the function $f'(z)/f(z)$ is almost periodic on each line $y=y_0\neq0$.
\end{Pro}

{\bf Proof}. Set
$$
\xi_z(t)=\begin{cases}-2\pi ie^{2\pi itz} &\text{if }t>0,\\0&\text{if }t\le0,\end{cases}\quad z\in\C_+,\quad\qquad
\xi_z(t)=\begin{cases}2\pi ie^{2\pi itz} &\text{if }t<0,\\0&\text{if }t\ge0,\end{cases}\quad z\in\C_-.
$$
It is not hard to check that in the sense of distributions $\hat\xi_z(\l)=1/(z-\l)$ for $z\in\C_+\cup\C_-$.
\medskip

Let $\p(t)$ be any  even nonnegative $C^\infty$-function such that $\supp\p\subset(-1,1)$ and $\int\p(t)dt=1$. Set $\p_\e(t)=\e^{-1}\p(t/\e)$
for $\e>0$. Fix $z=x+iy\in\C_+$. The functions $\xi_z(t)\star\p_\e(t)$ and $\hat\xi_z(\l)\hat\p_\e(\l)$ belong to $S(\R)$. Therefore,
$$
   (\hat\mu,\xi_z\star\p_\e(t))=(\mu,\hat\xi_z(\l)\hat\p_\e(\l)).
$$
Then for any  $T_0<\infty$ we have
$$
  \frac{i}{2\pi}(\hat\mu(t),\xi_z\star\p_\e(t))=d\int_{-\e}^0 e^{-2\pi isz}\p_\e(s)ds+ \sum_{0<|\g|\le\e}b_\g e^{2\pi i\g z}\int_{-\e}^\g e^{-2\pi isz}\p_\e(s)ds
$$
$$
 +\sum_{\e<\g<T_0}b_\g e^{2\pi i\g z}\int_{-\e}^\e e^{-2\pi isz}\p_\e(s)ds+\sum_{\g\ge T_0}b_\g e^{2\pi i\g z}\int_{-\e}^\e e^{-2\pi isz}\p_\e(s)ds=I_0+I_1+I_2+I_3.
$$
Then we have
$$
 \frac{i}{2\pi}(\hat\mu,\xi_z)=\sum_{0<\g\le\e}b_\g e^{2\pi i\g z}+\sum_{\e<\g<T_0}b_\g e^{2\pi i\g z}+\sum_{\g\ge T_0}b_\g e^{2\pi i\g z}=S_1+S_2+S_3.
$$
Set $m(s)=\sum_{\g\in\G:\,0<\g\le s}|b_\g|$. Then
\begin{equation}\label{s1}
   \sum_{\g\ge r}|b_\g|e^{-2\pi\g y}=\int_r^\infty e^{-2\pi sy}m(ds)\le\lim_{T\to\infty}m(T)e^{-2\pi Ty}+2\pi y\int_r^\infty e^{-2\pi sy}m(s)ds.
\end{equation}
Property \eqref{ch} implies $I_3$ and $S_3$ are less than given $\eta>0$ for $T_0$ large enough. Taking into account that $\sum_{-\e<\g<T_0}|b_\g|<\infty$, we get as $\e\to0$
$$
I_0\to d/2,\quad I_1\to0,\quad S_1\to 0, \quad  I_2-S_2=\sum_{\e<\g<T_0}b_\g e^{2\pi i\g z}\int_{-\e}^\e(e^{-2\pi isz}-1)\p_\e(s)ds\to0
$$
and
 $$
 \lim_{\e\to 0}(\mu,\hat\xi_z(\l)\hat\p_\e(\l))=(\hat\mu,\xi_z(t))-2\pi id/2=-2\pi i\sum_{\g\in\G\cap\R_+}b_\g e^{2\pi i\g z}-\pi id.
$$
 On the other hand, we have
\begin{equation}\label{h}
 (\mu,\hat\xi_z(\l)\hat\p_\e(\l))=\frac{\hat\p(\e a_0)}{z-a_0}+\sum_{n\in\N}\left[\frac{\hat\p(\e a_n)}{z-a_n}+\frac{\hat\p(\e a_{-n})}{z-a_{-n}}\right].
\end{equation}
The function $\hat\p(t)$ tends to $1$ as $t\to0$ and $|\hat\p(t)|\le1$.  We have
$$
 \left[\frac{\hat\p(\e a_n)}{z-a_n}+\frac{\hat\p(\e a_{-n})}{z-a_{-n}}\right]=\hat\p(\e a_{-n})\left[\frac{1}{z-a_n}+\frac{1}{z-a_{-n}}\right]+
 \frac{1}{z-a_n}[\hat\p(\e a_n)-\hat\p(\e a_{-n})].
$$
 Since $\hat\p$ is even, we get with bounded $\theta(n)$ and $\phi(n)$
$$
\hat\p(\e a_n)-\hat\p(\e a_{-n})=\hat\p(\e n/d+\e\phi(n))-\hat\p(-\e n/d+\e\phi(-n))=\hat\p'(\e n/d+\e\theta(n))\e[\phi(n)+\phi(-n)].
$$
Since $\hat\p(t)$ belongs to Schwartz space,  we see that $\hat\p'(t)=O(1/|t|)$ as $t\to\infty$. Hence for $\e|n/d+\theta(n)|>1$
$$
|\e[\hat\p'(\e n/d+\e\theta(n))]|\le C|n|^{-1}
$$
with a constant $C<\infty$. The same estimate (with another constant $C$) is valid for $\e|n/d+\theta(n)|\le1$,
that for all $n\in\N$ and $\e>0$
$$
|\hat\p(\e a_n)-\hat\p(\e a_{-n})|\le (C/n)2\sup_n|\phi(n)|.
$$
 Hence the right-hand side of \eqref{h} for all $\e>0$ is majorized by the  sum
\begin{equation}\label{sum}
\frac{1}{|z-a_0|}+\sum_{n\in\N}\left|\frac{1}{z-a_n}+\frac{1}{z-a_{-n}}\right|+\sum_{n\in\N}\frac{C'}{n|z-a_n|}.
\end{equation}
By Theorem \ref{T1}, we have  $1/(z-a_n)=O(1/n)$. Taking into account also Corollary \ref{C1}, we get the convergence of both sums in \eqref{sum}.
Therefore we can pass to the limit  in \eqref{h} as $\e\to0$ and obtain \eqref{n1}.

 By \eqref{s1}, $\sum_{\g\ge 1}|b_\g|e^{-2\pi\g y_0}<\infty$ for $y_0>0$,
and $\sum_{0<\g<1}|b_\g|<\infty$. Therefore the series in right-hand part of \eqref{n1} absolutely and uniformly for $\Im z\ge\a>0$ converges, and $f'(z)/f(z)$ is almost periodic on the line $y=y_0$.

In the case $y_0<0$ we apply \eqref{n1} to the function $\overline{f(\bar z)}$ and obtain \eqref{n2}. \bs

\begin{Th}\label{T3}  If
\begin{equation}\label{t3}
\int_0^1 s^{-2}|\hat\mu|(0,s)ds<\infty,
 \end{equation}
 then the function $f$ from \eqref{f} with the zero set $A$ can be rewritten in the form
$$
f(z)=\sum_{\o\in\O}q_\o e^{2\pi i\o z}\qquad q_\o\in\C\setminus\{0\}, \quad \o\in\R,\quad \sum_{\o\in\O}|q_\o|<\infty,
$$
where a countable bounded set $\O$ satisfies the conditions $\sup\O\in\O$, $\inf\O\in\O$.
\end{Th}
{\bf Proof}.
 The sum in the right-hand side of \eqref{n1} converges absolutely and uniformly in $x\in\R$ and $y\ge\a>0$. Changing the order of summing and integrating,  we get for $z=x+iy\in\C_+$
\begin{equation}\label{pr}
  \log f(z)-\log f(i)=\int_i^z\frac{f'(z)}{f(z)}dz=-\sum_{\g\in\G\cap\R_+}b_\g \frac{e^{2\pi i\g z}-e^{-2\pi\g}}{\g}-id\pi z-\pi d.
\end{equation}
It is easy to check that the convergence of the integral in \eqref{t3} implies (in fact is equivalent to) the convergence of the series
$$
\sum_{0<\g<1}|b_\g|\g^{-1}=\int_0^1 s^{-1}d|\hat\mu|(0,s),
 $$
and \eqref{s1} implies
$$
\sum_{\g\ge1}|b_\g|\g^{-1}e^{-2\pi\g}<\infty.
$$
Therefore,
\begin{equation}\label{pf}
  \log f(x+i)+id\pi x=-\sum_{0<\g<1}(b_\g/\g)e^{-2\pi\g}e^{2\pi i\g x}-\sum_{\g\ge1}(b_\g/\g)e^{-2\pi\g}e^{2\pi i\g x}+C_0
\end{equation}
with some constant $C_0\in\C$, and $\|\log f(x+i)+id\pi x\|_W<\infty$. Since $\|FG\|_W\le\|F\|_W\|G\|_W$,
we obtain
$$
   f(x+i)e^{id\pi x}=\sum_{\o\in\O}p_\o e^{2\pi i\o x},\quad \sum_{\o\in\O}|p_\o|=\|f(x+i)e^{id\pi x}\|_W\le e^{\|\log f(x+i)+id\pi x\|_W}<\infty,
 $$
with $p_\o\in\C$ and a countable spectrum $\O\subset\R_+\cup\{0\}$. The entire function $f(z+i)e^{id\pi z}$ has exponential growth, therefore, by \S 1, Ch.VI from \cite{L}, $\O$ is bounded.
Hence the function
$$
    f(z)=\sum_{\o\in\O}p_\o e^{\pi(2\o-d)}e^{\pi i(2\o-d)z}
$$
is also the Dirichlet series and $\|f\|_W<\infty$. Moreover, all zeros of $f$ are real, therefore, by Cor.2, Ch.VI from \cite{L}, we have $\sup\O\in\O$, $\inf\O\in\O$.  \bs

\medskip
Set
$$
g(z):=\sum_{\g\in\G,0<\g<1}b_\g\frac{e^{2\pi i\g z}-1}{\g}.
$$
 Since the sum $\sum_{\g\in\G,0<\g<1}|b_\g|$ is bounded and
$$
  g(z)=\sum_{k=1}^\infty\frac{(2\pi iz)^k}{k!}\sum_{\g\in\G,0<\g<1}\g^{k-1}b_\g,
$$
we see that $g(z)$ is a well-defined entire function.   Then for $z\in\C_+\cup\R$
\begin{equation}\label{g}
   |g(z)|\le\left[\sum_{\g\in\G,0<\g<\e}+\sum_{\g\in\G,\e\le\g<1}\right]\left|\frac{1- e^{2\pi i\g z}}{\g z}\right||z||b_\g|
\le4\pi|z|\sum_{0<\g<\e}|b_\g|+\frac{2}{\e}\sum_{\e\le\g<1}|b_\g|.
\end{equation}
The sum $\sum_{0<|\g|<\e}|b_\g|$ is arbitrary small for small $\e$, therefore, $|g(z)|=o(|z|)$ as $|z|\to\infty,\,z\in\C_+\cup\R$. Besides, in every strip $\{z:\,|\Im z|<M\}$
\begin{equation}\label{g1}
  |g(z)-g(x)|\le\sum_{\g\in\G,0<\g<1}|b_\g|\left|e^{2\pi i\g x}\frac{e^{-2\pi\g y}-1}{\g}\right|<C(M).
\end{equation}

\begin{Th}\label{T4}
If $g(z)$ is uniformly bounded for $z=x\in\R$, then the function \eqref{f} with the zero set $A$ is almost periodic of the exponential type $\pi d$.

Conversely, if $A$ is the zero set of an entire almost periodic function  of exponential growth, then the function $g(z)$ is uniformly bounded on $\R$.
\end{Th}

{\bf Proof}.
By \eqref{s1}, the sum $\sum_{\g\ge1}b_\g\g^{-1} e^{-2\pi\g}$ is finite, as is the sum $\sum_{0<\g<1}b_\g\frac{1-e^{-2\pi\g}}{\g}$. Hence
we can rewrite \eqref{pr} in the form
\begin{equation}\label{ph}
  \log f(z)=-\sum_{0<\g<1}b_\g \frac{e^{2\pi i\g z}-1}{\g}-\sum_{\g\ge1}b_\g\frac{e^{2\pi i\g z}}{\g}-id\pi z +\const.
\end{equation}
Here the second sum is uniformly bounded in every half-plane $\Im z\ge\a>0$. If the function $g(x)$ is uniformly bounded on $\R$, then by \eqref{g1},
the first sum in \eqref{ph} is bounded on every line $\Im z=M>0$. Therefore the functions $\log|f(z)|$ and $f(z)$ are bounded on this line, and $f(z)$ is bounded on every line $\Im z=-M<0$ as well.
Since the function $f(z)$ has the exponential growth, the Fragment-Lindelof Theorem implies boundedness of $f$ in every horizontal strip of bounded width.

Furthermore, combining \eqref{ph} and \eqref{g}, we get
$$
\lim_{y\to+\infty}y^{-1}\log|f(iy)|=\pi d.
$$
 The function $f(z)$ is bounded on the real line, therefore the Fragment-Lindelof Theorem implies that $|f(z)|\le Me^{y\pi d}$ for all $z\in\C_+$. Also, $|f(z)|\le Me^{-y\pi d}$ for all $z\in\C_-$,
 so $f(z)$ has the exponential type $\pi d$.

Next, the function $\log f(z)+id\pi z$ is bounded on the line $\Im z=1$. Its derivative  $(\log f(x+i))'+id\pi$ is almost periodic, hence, by Bohr's Theorem (cf.\cite{Le}, Theorem 1.2.1),
the functions $\log f(x+i)+id\pi x$ and $f(x+i)$ are almost periodic in the variable $x$.  By  Theorem 1.2.3 from \cite{Le}, the function $f(z)$ is almost periodic in every strip, where it is bounded,
hence it is almost periodic in $\C$.

Now let  $G(z)$ be an entire almost periodic function of exponential type with zero set $A$. Clearly, $G(z)=K_1e^{K_2z}f(z)$ with $K_1,\,K_2\in\C$.
  Taking into account that the second sum in \eqref{ph} is bounded on the line $\Im z=1$, we obtain as $x\to\infty$
 $$
 \log G(x+i)=K_2x-g(x+i)-id\pi x+O(1).
 $$
 Since $G(x+i)$ is almost periodic, left-hand side of this equality is also uniformly bound. Then \eqref{g1} implies
 $$
   K_2x-g(x)-id\pi x=O(1).
 $$
 By \eqref{g}, $g(x)=o(|x|)$, therefore, $K_2=i\pi d$ and $g(x)$ is bounded.  \bs

\bigskip
I would like to thank the Department of Mathematics and Computer Science
of the Jagiellonian University for their hospitality and Professor Lukasz Kosinski for his interest in my work.


\begin{thebibliography}{8}


\bibitem{B} H.Bohr, Almost Periodic Functions, ed. Chelsea,
New-York, 1951.

\bibitem{D} Directions in Mathematical Quasicrystals, M.Baake, R.Moody, eds.
CRM Monograph series 2000 {\bf 13}, AMS, Providence RI, 379p.

\bibitem{F} S.Yu.Favorov. Zeros of holomorphic almost periodic functions //Journal
d'Analyse Mathematique, v.84, (2001), p.51-66.

\bibitem{F1} S.Yu.Favorov. Uniqueness Theorems for Fourier Quasicrystals and Temperate Distributions with Discrete Support.
 Proc. Amer. Math. Soc. 149 (2021), 4431-4440

\bibitem{F4} S.Yu.Favorov {\it Large Fourier Quasicrystals and Wiener's Theorem},  Journal of Fourier Analysis and Applications, {\bf 25}, Issue 2, (2019), 377-392,
 DOI 10.1007/s00041-017-9576-0

 \bibitem{FK1} S.Yu.Favorov, Ye.Yu.Kolbasina. {\it Perturbations of discrete lattices and almost
periodic sets}. Algebra and Discrete Mathematica, {\bf 9} (2) (2010), 48-58.

\bibitem{FK2} S.Yu.Favorov, Ye.Yu.Kolbasina. Almost periodic discrete sets// Journal of
Mathematical Physics, Analysis, Geometry. v.6 (2010), No.1, 1-14.

\bibitem{FRR}S.Yu.Favorov, A.Yu.Rashkovskii, and L.I.Ronkin. Almost
periodic divisors in a strip //Journal d'Analyse Mathematique, Vol.74 (1998), 325-345.

\bibitem{La} W.Lawton {\it Bohr Almost Periodic Sets of Toral Type}. The Journal of Geometric Analysis 32:60, (2022).

\bibitem{L} B.Ja.Levin, Distributions of Zeros of Entire Functions. Transl. of Math. Monograph, Vol.5, AMS Providence, R1, 1980.

\bibitem{Le} B.M.Levitan, Almost periodic functions. Gostehizdat, 1953, Moskow. 396 p. (In Russian)

\bibitem{L1} J.C.Lagarias {\it Mathematical Quasicrystals and the Problem of Diffraction}, in \cite{D}, 61-93.

\bibitem{M1} Y.Meyer  {\it Quasicrystals, Almost Periodic Patterns, Mean--periodic Functions, and Irregular Sampling}, African Diaspora Journal of Mathematics,
 {\bf 13} no.1, (2012) 1-45.

\bibitem{M} Y.Meyer, {\it  Global and local estimates
on trigonometric sums}, Trans. R. Norw. Soc. Sci. Lett. 2018(2) 1-25.

\bibitem{OU} Olevskii, A., Ulanovskii A.  Fourier quasicrystals with unit masses// Comptes Rendus Mathématique, 2020, 358, no 11-12, p. 1207-1211
https://doi.org/10.5802/crmath.142

\bibitem{OU1} Olevskii, A., Ulanovskii A. A Simple Crystalline Measure.
arXiv:2006.12037v2, (2020).

\bibitem{Q} Quasicrystals and Discrete Geometry. J.Patera,ed., Fields Institute Monographs 1998, AMS, Providence RI, 289p.

\bibitem{R} L.I.Ronkin, L.I. {\it Almost Periodic Distributions and Divisors in Tube
Domains,} Zap. Nauchn. Sem. POMI {\bf 247} (1997)  210--236 (Russian).

\bibitem{T} Tornehave,H. {\it On the zeros of entire almost periodic function}.
The Harald Bohr Centenary (Copenhagen 1987). Math. Fys. Medd. Danske, {\bf 42}, no.3 (1989), 125-142.




\end{thebibliography}
\end{document}